\documentclass[10pt,psamsfonts]{amsart}
\usepackage{amsmath}
\usepackage{amsthm}
\usepackage{amssymb}
\usepackage{amscd}
\usepackage{amsfonts}
\usepackage{amsbsy}
\usepackage{graphicx}
\usepackage{url}
\usepackage{overpic}
\usepackage{hyperref}
\usepackage{float}
\usepackage{xcolor}
\usepackage[pagewise]{lineno}
\usepackage[top=2.5cm,bottom=2.5cm,left=2.5cm,right=2.5cm]{geometry}
\usepackage{enumitem}

\def\p{\partial}

\newtheorem {example} {Example}

\usepackage{mathpazo}

\begin{document}

\title[A note on a recent attempt to solve the second part of Hilbert's 16th Problem]
{A note on a recent attempt to solve\\ the second part of Hilbert's 16th Problem}

\author[C.A. Buzzi and D.D. Novaes]
{Claudio A. Buzzi$^1$ and Douglas D. Novaes$^2$}

\address{$^1$Universidade Estadual Paulista, IBILCE-UNESP - Av. Cristov\~ao Colombo, 2265, 15.054-000, S. J. Rio Preto, SP, Brasil}
\email{claudio.buzzi@unesp.br}
\address{$^2$Universidade Estadual de Campinas (UNICAMP), Departamento de Matem\'{a}tica, Instituto de Matemática, Estatística e Computação Científica (IMECC) - Rua S\'{e}rgio Buarque de Holanda, 651, Cidade Universit\'{a}ria Zeferino Vaz, 13083--859, Campinas, SP, Brasil}
\email{ddnovaes@unicamp.br}



\begin{abstract}
For a given natural number $n$, the second part of Hilbert’s 16th Problem asks whether there exists a finite upper bound for the maximum number of limit cycles that planar polynomial vector fields of degree $n$ can have. This maximum number of limit cycle, denoted by $H(n)$, is called the $n$th Hilbert number. It is well-established that $H(n)$ grows asymptotically as fast as $n^2 \log n$. A direct consequence of this growth estimation is that $H(n)$ cannot be bounded from above by any quadratic polynomial function of $n$. Recently, the authors of the paper [Exploring limit cycles of differential equations through information geometry unveils the solution to Hilbert’s 16th problem. Entropy, 26(9), 2024] affirmed to have solved the second part of Hilbert's 16th Problem by claiming that $H(n) = 2(n - 1)(4(n - 1) - 2)$. Since this expression is quadratic in $n$, it contradicts the established asymptotic behavior and, therefore, cannot hold. In this note, we further explore this issue by discussing some counterexamples.
\end{abstract}

\maketitle

\section{Introduction}

For a given natural number $n$, the second part of Hilbert's 16th Problem asks whether there is a finite upper bound for the number of limit cycles that planar polynomial vector fields of degree $n$ can possess. More precisely, let
\[
H(n) := \sup\{\pi(P, Q) : \deg(P), \deg(Q) \leq n\},
\]
where $\pi(P, Q)$ denotes the number of limit cycles of the polynomial differential system
\begin{equation}\label{ppvf}
\begin{cases}
\dot{x} = P(x, y), \\
\dot{y} = Q(x, y).
\end{cases}
\end{equation}
Recall that a limit cycle of \eqref{ppvf} is a (non-stationary) periodic orbit that is isolated from other periodic orbits (see \cite[Definition 9]{roussarie}). Thus, the second part of Hilbert's 16th Problem consists of proving that $H(n) < \infty$ for all $n \in \mathbb{N}$ (see \cite[Chapter 2]{roussarie}). The value $H(n)$ is called the $n$th Hilbert number.

The most significant advancement in understanding the asymptotic behavior of the function $H(n)$ was made by Christopher and Lloyd in \cite{CL95}, who introduced a method showing that $H(n)$ grows as fast as $n^2\log n$. This classical result has been revisited and improved by several works, including \cite{lvarez2020,HL12,LCC02}. In particular, Han and Li in \cite{HL12} refined Christopher and Lloyd's result, demonstrating that $H(n)$ grows at least as fast as $(n+2)^2\log (n+2)/(2\log 2)$ by establishing that
\[
\lim_{n\to\infty}\inf\dfrac{H(n)}{(n+2)^2\log(n+2)}\geq\dfrac{1}{2\log 2}.
\]
This remains the best-known lower estimation for the asymptotic growth  of $H(n)$.

A direct conclusion from this asymptotic growth  estimation is that $H(n)$ cannot be bounded from above by any quadratic polynomial function in $n$, as the expression $(n+2)^2\log (n+2)/(2\log 2)$ surpasses any degree two polynomial in $n$ for sufficiently large values of $n$.

Recently, the authors of the paper \cite{e26090745} affirmed to have solved the second part of Hilbert's 16th Problem by claiming that
\begin{equation}\label{upper}
H(n)= 2 (n - 1) (4 (n - 1) - 2),
\end{equation}
for $n\geq 2$ (see \cite[Theorem 4]{e26090745}). They make use of the following scalar curvature associated to a Fisher information metric:
\begin{equation}\label{R}
R=\dfrac{1}{\sqrt{G}}\left[ \dfrac{\p}{\p x}\left(\dfrac{1}{\sqrt{G}}\dfrac{\p G_{22}}{\p x}\right)+\dfrac{\p}{\p y}\left(\dfrac{1}{\sqrt{G}}\dfrac{\p G_{11}}{\p y}\right)\right],
\end{equation}
where 
\[
G_{11}=2\left[\left(\dfrac{\p P}{\p x}\right)^2+\left(\dfrac{\p Q}{\p x}\right)^2
\right],\,\, G_{22}=2\left[\left(\dfrac{\p P}{\p y}\right)^2+\left(\dfrac{\p Q}{\p y}\right)^2
\right],\,\,\text{and}\,\, G=G_{11} G_{22}.
\]
Their approach relies on \cite[Definition 1]{e26090745}, which aims to provide an alternative definition for limit cycles, referred to as being ``in the framework of GBT''. It begins by establishing  that
\begin{enumerate}[label = {$(A)$}]
\item\label{assertA} 
 {\it a limit cycle is the periodic state of \eqref{ppvf} in which $R$ is positive in the neighborhood of the equilibrium points of \eqref{ppvf} and $|R|$ is singular.}
\end{enumerate}
By $|R|$ singular, they mean the existence of zeros of the denominator of $|R|$ that makes $|R|$ to diverge to infinity. Thus, it is also asserted that  
\begin{enumerate}[label = {$(B)$}]
\item\label{assertB} 
{\it 
if $R$ is positive in the neighborhood of the equilibrium points of \eqref{ppvf} and the magnitude of $R$ diverges to infinity at symmetrical singularities with respect to the origin, then \eqref{ppvf} possesses only one limit cycle. Nonetheless, if $R$ is positive in the neighborhood of the equilibrium points of \eqref{ppvf} and the magnitude of $R$ diverges to infinity at different singularities, then \eqref{ppvf} has more than one limit cycle such that the total number of distinctive divergences of $|R|$ to infinity provides the maximum number of limit cycles of \eqref{ppvf}.}
 \end{enumerate}
Subsequent to Definition 1,  it is stated that  such a definition ``agrees with the definition of limit cycles in the framework of classical bifurcation theory'', that is (non-stationary) periodic orbits isolated from other periodic orbits. In this way, the approach employed in \cite{e26090745} to obtain \eqref{upper} consisted in counting the number of divergences of $|R|$ to infinity, as highlighted in the proof of \cite[Theorem 4]{e26090745}. 

As previously mentioned, the function $H(n)$ cannot be bounded from above by any quadratic polynomial in $n$. Therefore, the relationship \eqref{upper}, which is quadratic in $n$, cannot hold. To explore this issue further, we present counterexamples in the following sections. Section \ref{sec:kcl} discusses a well-known example from the literature that contradicts \eqref{upper}, along with references to other known examples that serves as counterexamples to \eqref{upper}. In Section \ref{sec:pipm}, we provide examples of polynomial systems that exhibit limit cycles but do not satisfy \ref{assertA}, and vice versa. This demonstrates that \ref{assertA} is neither necessary nor sufficient for the existence of limit cycles of \eqref{ppvf} and, therefore, is not equivalent to the standard definition of limit cycles. As a result, the definition of limit cycles proposed in \cite{e26090745} is not applicable to the study of the second part of Hilbert's 16th problem, meaning that the number of singularities of $|R|$ does not determine the maximum number of limit cycles in \eqref{ppvf}, as suggested by assertion \ref{assertB}.

\section{Known counterexamples in the literature}\label{sec:kcl}

The objective of this section is not to construct new counterexamples to the main conclusion \eqref{upper} of \cite{e26090745}, but rather to highlight known examples from the literature that serve as counterexamples for it.

In \cite[Section 3]{LCC02}, Li et al. revisited the class of polynomial differential systems originally studied by Christopher and Lloyd \cite{CL95}, addressing a minor issue in the original analysis. This correction did not affect the leading term $n^2\log n$ of the lower estimation for the asymptotic growth  of $H(n)$. Their approach, as well as Christopher and Lloyd's approach, consists of constructing a sequence of recursively defined polynomial differential systems $(PH_k)$ of degree $2^k - 1$, each possessing at least $S_k$ limit cycles, where
\[
S_k=4^{k-1}\left(k-\dfrac{13}{6}\right)+2^k-\dfrac{1}{3}.
\]
This sequence implies that
\begin{equation}\label{conclusion1}
H(2^k-1)\geq S_k=4^{k-1}\left(k-\dfrac{13}{6}\right)+2^k-\dfrac{1}{3}.
\end{equation}
However, the conclusion \eqref{upper} from \cite{e26090745} provides that
\[
H(2^k-1)=4 (2^k-2 ) (2^{k+1}-5 ),
\]
which contradicts \eqref{conclusion1} for $k\geq 35$. This means that system $PH_k$, for $k\geq 35$, has more limit cycles than predicted by the main result of \cite{e26090745}. The other sequences of polynomial systems discussed in \cite[Sections 4 and 5]{LCC02} also provide counterexamples to \eqref{upper}.

The works \cite{HL12} and, more recently, \cite{lvarez2020} also provide similar lower estimations for the asymptotic growth  of $H(n)$. Both works present sequences of polynomial differential systems with specified degrees and numbers of limit cycles, differing in the mechanisms used to generate these limit cycles. Counterexamples to \eqref{upper} can be derived from these sequences in a way analogous to the approach outlined above.

\section{Possible issue for the proposed method}\label{sec:pipm}

We begin by presenting three examples of polynomial differential systems where the existence of limit cycles is guaranteed, but assertion \ref{assertA} does not hold. Specifically, in these examples, either $R$ is negative in a neighborhood of the unique equilibrium point, or $R$ is positive in a neighborhood of the unique equilibrium point, but $|R|$ is not singular. These examples demonstrate that limit cycles satisfying \ref{assertA} do not encompass all possible limit cycles in polynomial systems. As a result, the maximum number of limit cycles satisfying \ref{assertA} for a polynomial system of degree $n$ does not provide an upper bound for $H(n)$. This likely explains why the main result \eqref{upper} of \cite{e26090745} does not agree with the established lower estimations for the asymptotic growth of $H(n)$, as discussed in the previous section.

\begin{example}\label{ex1}
We start by considering the following cubic vector field
\begin{equation}\label{s1}
\begin{cases}
\displaystyle \dot x= -y+x(x^2+y^2-1),\\
\displaystyle \dot y= x+y(x^2+y^2-1),
\end{cases}
\end{equation}
which has a single equilibrium point, located at the origin $(0,0)$. This vector field also has a unique limit cycle surrounding the origin. To see that, it is enough to write system \eqref{s1} in polar coordinates $(x,y)=(r\cos(\theta),r\sin(\theta))$ as follows:
\[
\begin{cases}
\displaystyle \dot r= r(r^2-1),\\
\displaystyle \dot \theta= 1.
\end{cases}
\]
This implies that system \eqref{s1} has a unique limit cycle which is unstable and whose orbit corresponds to the unit circle with center at the origin. Now, computing the function $R$ we get
\[
R(x,y)=\dfrac{R_1(x,y)}{R_2(x,y)},
\]
where
\[
\begin{aligned}
R_1(x,y)=&72 x^{10}-216 x^8 y^2-204 x^8-320 x^7 y-3056 x^6 y^4+464 x^6 y^2+368 x^6+192 x^5 y^3+192 x^5 y\\
&-3056 x^4 y^6+2360 x^4 y^4-304 x^4 y^2-240 x^4-192 x^3 y^5-216 x^2 y^8+464 x^2 y^6-304 x^2 y^4-96 x^2 y^2\\
&+96 x^2+320 x y^7-192 x y^5+72 y^{10}-204 y^8+368 y^6-240 y^4+96 y^2-16\quad \text{and}\quad\\
R_2(x,y)=& \Big((3 x^2+y^2-1)^2+(2 x y+1)^2\Big)^2 \Big((x^2+3 y^2-1)^2+(2 x y-1)^2\Big)^2.
\end{aligned}
\]
Observe that $R_2$ does not vanish at the origin, implying that $R$ is continuous in its neighborhood. Additionally, since $R(0,0) = -1 < 0$, continuity ensures that $R(x,y)$ remains negative in a neighborhood of the origin, which corresponds to the unique equilibrium point of \eqref{s1}. Therefore, system \eqref{s1} provides an example of a limit cycle that does not satisfy assertion \ref{assertA}.
\end{example}

\begin{example}\label{ex2}
Using the approach from Example \ref{ex1}, we can easily construct polynomial systems with any number of limit cycles and a unique equilibrium point, where $R$ is negative in its neighborhood. For instance, the following polynomial system has a single equilibrium point at the origin and two nested limit cycles surrounding it:
\begin{equation}\label{s1a}
\begin{cases}
\displaystyle \dot x= -y+x(x^2+y^2-1)(x^2+y^2-4),\\
\displaystyle \dot y= x+y(x^2+y^2-1)(x^2+y^2-4).
\end{cases}
\end{equation}
Indeed, by applying a polar change of variables, one can deduce that \eqref{s1a} has exactly two limit cycles: an asymptotically stable one, whose orbit corresponds to the unit circle centered at the origin; and an unstable one whose orbit corresponds  to a circle of radius two, also centered at the origin. The expression for $R$ is cumbersome and thus omitted here, but following the same reasoning of Example \ref{ex1}, we conclude that $R$ is continuous in a neighborhood of the origin, with $R(0,0) = -80/289 < 0$, implying that $R$ remains negative near the origin. Therefore, system \eqref{s1a} provides examples of limit cycles that do not satisfy assertion \ref{assertA}.
\end{example}

\begin{example}
Now, consider the system \eqref{s1} under the following linear change of variables: $(x, y) = (u, u + v/2)$. This yields the transformed system:
\begin{equation}\label{s2}
\begin{cases}
\displaystyle \dot u= -2 u-\frac{v}{2}+2 u^3+u^2 v+\frac{u v^2}{4},\\
\displaystyle \dot v=4 u+ 2 u^2 v+u v^2+\frac{v^3}{4}.
\end{cases}
\end{equation}
Of course, system \eqref{s2} has a unique equilibrium point at the origin $(0,0)$ and a unique limit cycle surrounding it. Computing the function $R$ for system \eqref{s2}, we obtain
\[
R(u,v)=\dfrac{R_1(u,v)}{R_2(u,v)},
\]
where
\[
\begin{aligned}
R_1(u,v)=&32\Big(-663552 u^{10}-8638464 u^9 v-25353216 u^8 v^2-7421952 u^8-37943808 u^7 v^3-18733056 u^7 v\\
&-36060032 u^6 v^4-22151168 u^6 v^2+5670912 u^6-23658048 u^5 v^5-18140416 u^5 v^3+10874880 u^5 v\\
&-10971920 u^4 v^6-11152128 u^4 v^4+7196416 u^4 v^2-2199552 u^4-3555048 u^3 v^7-4852576 u^3 v^5\\
&+2186496 u^3 v^3-4174848 u^3 v-772632 u^2 v^8-1359232 u^2 v^6+296160 u^2 v^4-2595840 u^2 v^2\\
&+219136 u^2-103056 u v^9-222052 u v^7+49248 u v^5-828032 u v^3+472064 u v-6399 v^{10}-18528 v^8\\
&+18596 v^6-126272 v^4+134912 v^2+61440\Big) \quad \text{and}\quad\\
R_2(u,v)=& \Big((24 u^2+8 u v+v^2-8)^2+16 (4 u v+v^2+4)^2\Big)^2 \Big((8 u^2+8 u v+3 v^2)^2+4 (2 u^2+u v-1)^2\Big)^2.
\end{aligned}
\]
Again, $R_2$ does not vanish at the origin, so $R$ is continuous in its neighborhood. Moreover, since $R(0,0) = 6/5 > 0$, continuity ensures that $R(u,v)$ is positive in a neighborhood of the origin, corresponding to the unique equilibrium point of \eqref{s2}. Additionally, since $R_2$ is a product of sums of squares, it follows that $R_2(u,v) = 0$ if and only if $(u,v)$ satisfies one of the following systems of algebraic equations:
\[
S_1:\begin{cases}
24 u^2+8 u v+v^2-8=0\\
4 u v+v^2+4=0
\end{cases}\quad \text{or}\quad S_2:\begin{cases}
8 u^2+8 u v+3 v^2=0\\
2 u^2+u v-1=0.
\end{cases}
\]
We begin by analyzing $S_1$. First, note that if $(u,v)$ is a solution of $S_1$, then $v \neq 0$. Solving the second equation of $S_1$ for $u$ and substituting into the first equation yields the algebraic equation $17 v^4 + 152 v^2 + 384 = 0$, which has no real solutions. Next, for system $S_2$, if $(u,v)$ is a solution, then $u \neq 0$. Solving the second equation of $S_2$ for $v$ and substituting into the first equation leads to the algebraic equation $3 - 4 u^2 + 4 u^4 = 0$, which also has no real solutions. This shows that the denominator $R_2$ of $R$ does not vanish, and hence $|R|$ has no singularities. Therefore, system \eqref{s2} provides another example of a limit cycle that does not satisfy assertion \ref{assertA}.
\end{example}

From the above examples, we observed that assertion \ref{assertA} is not necessary for the existence of limit cycles, as there are polynomial systems with limit cycles where \ref{assertA} does not hold. Nevertheless, we can still ask whether \ref{assertA} is a sufficient condition for the existence of limit cycles. The following example provides a negative answer to this question.

\begin{example}\label{ex-center}
Consider the following quadratic polynomial system:
\begin{equation}\label{center}
\begin{cases}
\displaystyle \dot x=-y+x^2,\\
\displaystyle \dot y= x+x y.
\end{cases}
\end{equation}
This system and its properties have been extensively studied in the literature, as it appears as a normal form for a class of isochronous quadratic systems, commonly referred to as $S_2$ (see \cite{loud,mardesic}). This system has a unique equilibrium point at the origin, which is a center, meaning that there exists a neighborhood $U$ around the origin where all orbits in $U \setminus \{(0,0)\}$ are periodic. Clearly, no periodic orbit in $U$ is a limit cycle, as none are isolated from other periodic orbits. In fact, this system does not have any limit cycles. By computing the function $R$ for system \eqref{center}, we obtain
\[
R(x,y)=\frac{1}{\left(x^2+1\right)^2 \left(4 x^2+(y+1)^2\right)}.
\]
Observe that $R$ is continuous in a neighborhood of the origin, as its denominator does not vanish at $(0,0)$. Since $R(0,0) = 1$, continuity ensures that $R$ remains positive in a neighborhood around the origin, which is the unique equilibrium point of \eqref{center}. Furthermore, $|R|$ is singular at $(x,y) = (0,-1)$. Thus, system \eqref{center} provides an example where assertion \ref{assertA} holds for every periodic orbit within $U$, despite the absence of limit cycles.
\end{example}

\section{Conclusion}

In this note, we have demonstrated that the recent attempt to solve the second part of Hilbert’s 16th problem, as presented in \cite{e26090745}, contains significant issues. We began by exploring counterexamples which demonstrate that the quadratic expression proposed for $H(n)$ contradicts the well-established asymptotic behavior of this function, which states that $H(n)$ grows as fast as $(n+2)^2 \log (n+2)/(2 \log 2)$. Moreover, we discussed how the alternative definition of limit cycles  \ref{assertA}, used in \cite{e26090745}, is not applicable to the study of the second part of Hilbert’s 16th problem, as it is neither necessary nor sufficient for the existence of limit cycles in \eqref{ppvf}, according to the standard definition, which refers to (non-stationary) periodic orbits isolated from other periodic orbits.

\bibliographystyle{abbrv}
\bibliography{references.bib}

\begin{thebibliography}{1}

\bibitem{lvarez2020}
M.~{\'{A}}lvarez, B.~Coll, P.~D. Maesschalck, and R.~Prohens.
\newblock Asymptotic lower bounds on {H}ilbert numbers using canard cycles.
\newblock {\em Journal of Differential Equations}, 268(7):3370--3391, Mar.
  2020.

\bibitem{CL95}
C.~J. Christopher and N.~G. Lloyd.
\newblock Polynomial systems: a lower bound for the {H}ilbert numbers.
\newblock {\em Proc. Roy. Soc. London Ser. A}, 450(1938):219--224, 1995.

\bibitem{e26090745}
V.~B. da~Silva, J.~P. Vieira, and E.~D. Leonel.
\newblock Exploring limit cycles of differential equations through information
  geometry unveils the solution to {H}ilbert’s 16th problem.
\newblock {\em Entropy}, 26(9), 2024.

\bibitem{HL12}
M.~Han and J.~Li.
\newblock Lower bounds for the {H}ilbert number of polynomial systems.
\newblock {\em J. Differential Equations}, 252(4):3278--3304, 2012.

\bibitem{LCC02}
J.~Li, H.~S.~Y. Chan, and K.~W. Chung.
\newblock Some lower bounds for {$H(n)$} in {H}ilbert's 16th problem.
\newblock {\em Qual. Theory Dyn. Syst.}, 3(2):345--360, 2002.

\bibitem{loud}
W.~S. Loud.
\newblock Behavior of the period of solutions of certain plane autonomous
  systems near centers.
\newblock {\em Contributions to Differential Equations}, 3:21--36, 1964.

\bibitem{mardesic}
P.~Marde\v{s}i\'{c}, C.~Rousseau, and B.~Toni.
\newblock Linearization of isochronous centers.
\newblock {\em J. Differential Equations}, 121(1):67--108, 1995.

\bibitem{roussarie}
R.~Roussarie.
\newblock {\em Bifurcation of planar vector fields and {H}ilbert's sixteenth
  problem}, volume 164 of {\em Progress in Mathematics}.
\newblock Birkh\"auser Verlag, Basel, 1998.

\end{thebibliography}
\end{document}